\documentclass[10pt]{article}
\usepackage{amssymb}

\newtheorem{prop}{Proposition}[section]
\newtheorem{thm}[prop]{Theorem}

\title{Simulation of the Collatz 3$x$\ +\ 1 function\\ by Turing machines}
\author{Pascal MICHEL\thanks{Corresponding address: 59 rue du Cardinal Lemoine,
75005 Paris, France.}\\
{\small \'Equipe de Logique Math\'ematique,}\\
{\small Institut de Math\'ematiques de Jussieu -- Paris Rive Gauche, UMR 7586,}\\
{\small B\^atiment Sophie Germain, case 7012, 75205 Paris Cedex 13, France}\\
{\small and Universit\'e de Cergy-Pontoise, IE, F-95000 Cergy-Pontoise, France}\\
{\small michel@math.univ-paris-diderot.fr}}

\date{September 24, 2014}

\begin{document}
\maketitle

\begin{abstract}
We give new Turing machines that simulate the iteration of the Collatz
$3x+1$ function. First, a never halting Turing machine with 3 states and
4 symbols, improving the known $3 \times 5$ and $4 \times 4$ Turing machines.
Second, Turing machines that halt on the final loop, in the classes $3 \times 10$,
$4 \times 6$, $5 \times 4$ and $13 \times 2$.

\bigskip


\bigskip

\end{abstract}

\section{Introduction}
Turing machines can be classified according to their numbers of states and symbols.
It is known (see \cite{WN09} for a survey) that there are universal Turing 
machines in the following sets (number of states $\times$ number of symbols):
$$2 \times 18,\ 3 \times 9,\ 4 \times 6,\ 5 \times 5,\ 6 \times 4,\ 9 \times 3,\ 18 \times 2.$$
On the other hand, all the Turing machines in the following sets are decidable:
$$1 \times n,\ 2 \times 3,\ 3 \times 2, n \times 1.$$
In order to refine the classification of Turing machines between universal and decidable classes,
properties in connection with the $3x + 1$ function have been considered.

Recall that the $3x + 1$ function $T$ is defined by
$$T(x) =\left\{
\begin{array}{ll}
x/2       & \mbox{ if $x$ is even}\\
(3x +1)/2 & \mbox{ if $x$ is odd}
\end{array}\right.$$
This can also be written $T(2n) = n$, $T(2n + 1) = 3n + 2$.
When function $T$ is iterated on a positive integer, it seems that the loop
$2 \mapsto 1 \mapsto 2$ is always reached, but this is unproven, and is a famous open problem
in mathematics \cite{La10}. For further references, we set
\begin{quote}
{\bf {\boldmath $3 x + 1$} Conjecture}: When function $T$ is iterated from positive integers,
the loop $2 \mapsto 1 \mapsto 2$ is always reached.
\end{quote}
The $3x + 1$ function is also called the Collatz function, and \emph{Collatz-like} functions
are functions on integers with a definition of the following form: 
there exist integers $d \ge 2$, $a_i$, $b_i$, $0 \le i \le d-1$, such that, for all integers $x$,
$$f(x) = \frac{a_ix + b_i}{d}\quad \mbox{if}\quad x \equiv i\quad \mbox{(mod $d$)}.$$
With these definitions, we can state the following properties of Turing machines,
that have been used to refine the classification according to the numbers of states and symbols
(see \cite{MM10} for a survey).
\begin{itemize}
\item Turing machines that simulate the iteration of the $3x + 1$ function
and never halt. It is known that there are such machines in the sets
$$2 \times 8,\ 3 \times 5,\ 4 \times 4,\ 5 \times 3,\ 10 \times 2.$$
We improve these results by giving a $3 \times 4$ Turing machine.
\item Turing machines that simulate the iteration of the $3x +1$ function
and halt when the loop $2 \mapsto 1 \mapsto 2$ is reached.
It is known that there is such a machine in the set $6 \times 3$.
In this article, we give four new Turing machines, in the classes
$3 \times 10$, $4 \times 6$, $5 \times 4$ and $13 \times 2$.
\item Turing machine that simulate the iteration of a Collatz-like function.
It is known that there are such machines in the sets
$$2 \times 4,\ 3 \times 3,\ 5 \times 2.$$
\end{itemize}

\section{Preliminaries: Turing machines}
The Turing machines we use have
\begin{itemize}
\item one tape, infinite on both sides, made of cells containing symbols,
\item one reading and writing head,
\item a set $Q = \{A, B, \ldots\}$ of states, plus a halting state $H$ (or $Z$),
\item a set $\Sigma = \{b,0,1,\ldots\}$ of symbols, where $b$ is the blank symbol
(or $\Sigma = \{0,1\}$, when 0 is the blank symbol),
\item a next move function
$$\delta : Q \times \Sigma \rightarrow \Sigma \times \{L, R\} \times(Q \cup\{H\}).$$
\end{itemize}
If $\delta(p,a) = (b,D,q)$, then the Turing machine, reading symbol $a$ in state $p$,
replaces $a$ by $b$, moves in the direction $D \in \{L, R\}$ ($L$ for Left, 
$R$ for Right), and comes into state $q$.
On an input $x_k\ldots x_0 \in \Sigma^{k+1}$, the initial configuration is
$^\omega b(Ax_k)\ldots x_0b^\omega$. This means that the word $x_k\ldots x_0$ is written
on the tape between two infinite strings of blank symbols, and the machine is reading
symbol $x_k$ in state $A$.

\begin{table}
$$\begin{array}{c|c|c|c|c|c|c|c|c|c|c|c|c|c}
\mbox{symbols} & \multicolumn{13}{c}{}\\                                                \cline{1-3}
10             & Ma    & \mbox{\bf Mi}_2& \multicolumn{11}{c}{}\\                       \cline{1-3}
9              &       &      &\multicolumn{11}{c}{}\\                                  \cline{1-3}
8              & Ba    &      &\multicolumn{11}{c}{}\\                                  \cline{1-3}
7              &       &      &\multicolumn{11}{c}{}\\                                  \cline{1-4}
6              &       & Ma   &  \mbox{\bf Mi}_2 & \multicolumn{10}{c}{}\\          \cline{1-4}
5              &       & Ba   &                  & \multicolumn{10}{c}{}\\          \cline{1-5}
4              &       & Mi_2 & Ma               &\mbox{\bf Mi}_2 & \multicolumn{9}{c}{}\\         \cline{1-6}
3              &       &      &                  & Ma   &\mbox{\bf Mi}_1 & \multicolumn{8}{c}{}\\  \cline{1-13}
2              &       &      &    &      &      &     &     &     & Ba  & Ma & & \mbox{\bf Mi}_2\\  \hline
               & 2     & 3    & 4  & 5    &\,6\, &\,7\,&\,8\,&\,9\,& 10  & 11 & 12 & 13 & \mbox{states}
\end{array}$$
\caption{Turing machines simulating the $3x + 1$ function:
$Ma=$ Margenstern \cite{Ma98,Ma00},
$Ba=$ Baiocchi \cite{Ba98},
$Mi_1=$ Michel \cite{Mi93},
$Mi_2=$ Michel (this paper).
In roman boldface, halting machines.}
\end{table}

\section{The known Turing machines}
Let us give some more precisions about the Turing machines that simulate
the $3x + 1$ function. The following results are displayed in Table 1.

Michel \cite{Mi93} gave a $6 \times 4$ Turing machine that halts when number 1 is reached.
This machine works on numbers written in binary. Division by 2 of even integers
is easy and multiplication by 3 is done by the usual multiplication algorithm.

Margenstern \cite{Ma98,Ma00} gave never halting $5 \times 3$ and $11 \times 2$
Turing machines in binary, and never halting $2 \times 10$, $3 \times 6$,
$4 \times 4$ Turing machines in unary, that is working on numbers $n$ written
as strings of $n$ 1s.

Baiocchi \cite{Ba98} gave five never halting Turing machines in unary, including
$2 \times 8$, $3 \times 5$ and $10 \times 2$ machines that improved
Margenstern's results.

In this article, we give a never halting $3 \times 4$ Turing machine
that works on numbers written in base 3. Multiplication by  3 is easy
and division by 2 is done by the usual division algorithm.
Note that Baiocchi and Margenstern \cite{BM01} already used numbers written in base 3
to define cellular automata that simulate the $3x + 1$ function.

By adding two states to this $3 \times 4$ Turing machine,
we derive a $5 \times 4$ Turing machine that halts when number 1 is reached.

We also give three other Turing machines that halt when number 1 is reached:

\begin{itemize}
\item A $3 \times 10$ Turing machine obtained by adding one state to the
$2 \times 10$ Turing machine of Margenstern \cite{Ma98,Ma00}.
\item A $4 \times 6$ Turing machine obtained by adding one state to the
$3 \times 6$ Turing machine of Margenstern \cite{Ma98,Ma00}.
\item A $13 \times 2$ Turing machine obtained by adding two states to the
$11 \times 2$ Turing machine of Margenstern \cite{Ma98,Ma00}.
\end{itemize}

\section{A never halting $3 \times 4$ Turing machine}
This Turing machine $M_1$ is defined as follows:

\begin{center}
\begin{tabular}{|c|c|c|c|c|}
\hline
$M_1$ &  $b$  &   0   &   1   &   2 \\
\hline
$A$ & $b$L$C$ & 0R$A$ & 0R$B$ & 1R$A$ \\
\hline
$B$ &  2L$C$  & 1R$B$ & 2R$A$ & 2R$B$ \\
\hline
$C$ & $b$R$A$ & 0L$C$ & 1L$C$ & 2L$C$ \\
\hline
\end{tabular}
\end{center}

The idea is simple. A positive integer is written on the tape, in base 3, in the usual order.
Initially, in state $A$, the head reads the most significant digit, at the left end
of the number. The initial configuration on input $x = \sum_{i=0}^k x_i3^i$ is
$^\omega b(Ax_k)\ldots x_0b^\omega$. Then the machine performs the division by 2, using the usual
division algorithm. Partial quotients are written on the tape. Partial remainders
are stored in the states: 0 in state $A$, 1 in state $B$. When the head passes the right
end of the number, reading a $b$, then
\begin{itemize}
\item if the remainder is 0, nothing is done: $2n \mapsto n$,
\item if the remainder is 1, a 2 is concatenated to the number:
$2n + 1 \mapsto n \mapsto 3n + 2$.
\end{itemize}
Then the head comes back, in state $C$, to the left end of the number and is
ready to perform a new division by 2.

We have the following theorem.
\begin{thm}
The $3x + 1$ conjecture is true iff, for all positive integer
$x = x_k\ldots x_0$ written in base 3, there exists an integer
$n \ge 0$ such that, on input $x_k\ldots x_0$, the Turing machine
$M_1$ eventually reaches the configuration ${^\omega}b0^n(A1)b^\omega$.
\end{thm}

\section{Turing machines that halts on the final loop}
\subsection{A $3 \times 10$ Turing machine}
Margenstern \cite[Fig.\ 11]{Ma00} gave the folowing never halting $2 \times 10$
Turing machine $M_2$.

\begin{center}
\begin{tabular}{|c|c|c|c|c|c|c|c|c|c|c|}
\hline
$M_2$ &  $b$  &    1    &   $x$   &   $r$   &   $u$   
    &  $v$  &    $y$    &   $z$   &   $t$   &   $k$   \\
\hline
$A$ & $b$R$A$ & $x$R$B$ &  1L$A$  & $k$R$B$ & $x$R$A$ 
& $x$R$A$ & $r$L$A$ & $r$L$A$ & $y$R$A$ &         \\
\hline
$B$ & $z$L$B$ & $u$R$B$ & $x$R$B$ & $y$R$B$ & $v$L$B$ 
& $u$R$A$ & $t$L$B$ &  1L$A$  & $x$R$B$ & $b$R$B$ \\
\hline
\end{tabular}
\end{center}

Turing machine $M_2$ works on numbers written in unary, so that the initial
configuration on number $n \ge 1$ is $^\omega b(A1)1^{n-1}b^\omega$.
By adding a new state $C$, we can detect the partial configuration
$(A1)b$, and we obtain the following $3 \times 10$ Turing machine $M_3$.

\begin{center}
\begin{tabular}{|c|c|c|c|c|c|c|c|c|c|c|}
\hline
$M_3$ &  $b$  &    1    &   $x$   &   $r$   &   $u$   
    &  $v$  &    $y$    &   $z$   &   $t$   &   $k$   \\
\hline
$A$ & $b$R$A$ & $x$R$C$ &  1L$A$  & $k$R$B$ & $x$R$A$ 
& $x$R$A$ & $r$L$A$ & $r$L$A$ & $y$R$A$ &         \\
\hline
$B$ & $z$L$B$ & $u$R$B$ & $x$R$B$ & $y$R$B$ & $v$L$B$ 
& $u$R$A$ & $t$L$B$ &  1L$A$  & $x$R$B$ & $b$R$B$ \\
\hline
$C$ & $b$L$H$ & $u$R$B$ &         & $y$R$B$ &         
 &         &         &         &         &         \\
\hline
\end{tabular}
\end{center}

We have the following theorem
\begin{thm}
The $3x + 1$ conjecture is true iff, for all positive integers $n$,
Turing machine $M_3$ halts on the initial configuration $^\omega b(A1)1^{n-1}b^\omega$.
\end{thm}

\subsection{A $4 \times 6$ Turing machine}
Margenstern \cite[Fig.\ 10]{Ma00} gave the folowing never halting $3 \times 6$
Turing machine $M_4$ (Note that transition $(1,z) \mapsto (x\mbox{R}2)$
in this figure should be $(1,z) \mapsto (r\mbox{R}2)$).

\begin{center}
\begin{tabular}{|c|c|c|c|c|c|c|}
\hline
$M_4$ &  $b$  &    1    & $x$   &    $a$  &   $z$   &  $r$ \\
\hline
$A$ & $b$R$A$ & $x$R$B$ & 1L$A$ &  1L$A$  & $r$R$B$ & \\
\hline
$B$ &  1L$B$  & $a$R$C$ & 1L$B$ &  1L$A$  & $x$R$B$ & $b$R$A$ \\
\hline
$C$ & $z$L$C$ & $x$R$C$ & 1L$C$ & $a$R$A$ & $r$R$C$ & $z$L$C$ \\
\hline
\end{tabular}
\end{center}

Turing machine $M_4$ works on numbers written in unary, with initial
configuration ${^\omega}b(A1)1^{n-1}b^\omega$.
By adding a new state $D$, we can detect
the partial configuration $(A1)b$, and we obtain the
following $4 \times 6$ Turing machine $M_5$.

\begin{center}
\begin{tabular}{|c|c|c|c|c|c|c|}
\hline
$M_5$ &  $b$  &    1    & $x$   &    $a$  &   $z$   &  $r$ \\
\hline
$A$ & $b$R$A$ & $x$R$D$ & 1L$A$ &  1L$A$  & $r$R$B$ & \\
\hline
$B$ &  1L$B$  & $a$R$C$ & 1L$B$ &  1L$A$  & $x$R$B$ & $b$R$A$ \\
\hline
$C$ & $z$L$C$ & $x$R$C$ & 1L$C$ & $a$R$A$ & $r$R$C$ & $z$L$C$ \\
\hline
$D$ & $b$L$H$ & $a$R$C$ &       &         & $x$R$B$ &         \\
\hline
\end{tabular}
\end{center}

We have the following theorem.
\begin{thm}
The $3x + 1$ conjecture is true iff, for all positive integers $n$,
Turing machine $M_5$ halts on the initial configuration ${^\omega}b(A1)1^{n-1}b^\omega$.
\end{thm}

\subsection{A $5 \times 4$ Turing machine}
This Turing machine $M_6$ is defined as follows.

\begin{center}
\begin{tabular}{|c|c|c|c|c|}
\hline
$M_6$ &  $b$  &    0    &    1    &   2 \\
\hline
$A$ & $b$L$C$ &  0R$A$  &  0R$B$  & 1R$A$ \\
\hline
$B$ &  2L$E$  &  1R$B$  &  2R$A$  & 2R$B$ \\
\hline
$C$ & $b$R$D$ &  0L$C$  &  1L$C$  & 2L$C$ \\
\hline
$D$ &         & $b$R$A$ & $b$R$B$ & 1R$A$ \\
\hline
$E$ & $b$R$H$ &  0L$C$  &  1L$C$  & 2L$C$ \\
\hline
\end{tabular}
\end{center}

Turing machine $M_6$ is obtained from Turing machine $M_1$ by adding
a state $D$ that wipes out the useless 0s, and a state $E$ that detects
the partial configuration $b(Bb)$.

We have the following theorem.
\begin{thm}
The $3x + 1$ conjecture is true iff Turing machine $M_6$ halts on all
input $x = x_k\ldots x_0$ representing a positive integer written in base 3.
\end{thm}

\subsection{A $13 \times 2$ Turing machine}
Margenstern \cite[Fig.\ 8]{Ma00} gave the following never halting $11 \times 2$
Turing machine $M_7$ (in this table, $H$ is \emph{not} a halting state).

\begin{center}
\begin{tabular}{|c|c|c|}
\hline
$M_7$ &  0  &   1   \\
\hline
$A$ & 1R$I$ & 0R$B$ \\
\hline
$B$ & 0R$A$ & 0R$G$ \\
\hline
$C$ & 0R$A$ & 1R$D$ \\
\hline
$D$ & 0R$C$ & 1R$E$ \\
\hline
$E$ & 1R$I$ & 1R$F$ \\
\hline
$F$ & 1R$C$ & 0R$G$ \\
\hline
$G$ & 1R$C$ & 1R$H$ \\
\hline
$H$ & 0R$E$ & 1R$G$ \\
\hline
$I$ & 1L$J$ &       \\
\hline
$J$ & 0R$B$ & 1L$K$ \\
\hline
$K$ & 0L$J$ & 1L$J$ \\
\hline
\end{tabular}
\end{center}

This machine works on numbers written in binary, with the least significant bit
at the left end of the number, and digits 0 and 1 coded by 10 and 11,
so that the initial configuration on number $n = x_k\ldots x_0 = \sum_{i=0}^k x_i2^i$
is $^\omega 0(A1)x_01x_1\ldots 1x_k0^\omega$. Division by 2 of even integers 
is easy, and multiplication by 3 is done by the usual algorithm.

By adding two new states $L$ and $M$, we can detect the partial configuration
$(A1)10$, and we obtain the following $13 \times 2$ Turing machine $M_8$,
where $Z$ is the halting state.

\begin{center}
\begin{tabular}{|c|c|c|}
\hline
$M_8$ &  0  &   1   \\
\hline
$A$ & 1R$I$ & 0R$L$ \\
\hline
$B$ & 0R$A$ & 0R$G$ \\
\hline
$C$ & 0R$A$ & 1R$D$ \\
\hline
$D$ & 0R$C$ & 1R$E$ \\
\hline
$E$ & 1R$I$ & 1R$F$ \\
\hline
$F$ & 1R$C$ & 0R$G$ \\
\hline
$G$ & 1R$C$ & 1R$H$ \\
\hline
$H$ & 0R$E$ & 1R$G$ \\
\hline
$I$ & 1L$J$ &       \\
\hline
$J$ & 0R$B$ & 1L$K$ \\
\hline
$K$ & 0L$J$ & 1L$J$ \\
\hline
$L$ & 0R$A$ & 0R$M$ \\
\hline
$M$ & 0L$Z$ & 1R$H$ \\
\hline
\end{tabular}
\end{center}

We have the following theorem.
\begin{thm}
The $3x + 1$ conjecture is true iff, for all positive number
$n = x_k\ldots x_0 = \sum_{i=0}^k x_i2^i$, Turing machine $M_8$ halts on the
initial configuration $^\omega 0(A1)x_01x_1\ldots 1x_k0^\omega$.
\end{thm}

\section{Conclusion}
We have given a new $3 \times 4$ never halting Turing machine
that simulates the iteration of the $3x + 1$ function.
It seems that it will be hard to improve the known results
on never halting machines.

On the other hand, for Turing machines that halt on the conjectured
final loop of the $3x + 1$ function, more researches are still to be done.


\begin{thebibliography}{99}
\bibitem{Ba98} C.\ Baiocchi, 3N+1, UTM e Tag-systems (Italian),
Dipartimento di Matematica dell'Universit\`a ``La Sapienza'' di Roma {\bf 98/38}, 1998.

\bibitem{BM01} C.\ Baiocchi and M.\ Margenstern, Cellular automata about
the $3x + 1$ problem, in: Proc.\ LCCS'2001, Universit\'e Paris 12, 2001, 37--45,
available on the website http://lacl.univ-paris12.fr/LCCS2001/.

\bibitem{La10} J.C.\ Lagarias (Ed.), The Ultimate Challenge: The 3$x$+1 Problem, AMS, 2010.

\bibitem{Ma98} M.\ Margenstern, Frontier between decidability and undecidability: a survey,
Proc.\ MCU'98, Vol.\ 1, ISBN 2-9511539-2-9, 1998, 141--177.

\bibitem{Ma00} M.\ Margenstern, Frontier between decidability and undecidability:
a survey, \emph{Theoret.\ Comput.\ Sci.} {\bf 231}, 2000, 217--251.

\bibitem{Mi93} P.\ Michel, Busy beaver competition and Collatz-like problems,
\emph{Arch.\ Math.\ Logic} {\bf 32} (5), 1993, 351--367.

\bibitem{MM10} P.\ Michel and M.\ Margenstern, Generalized 3$x$+1 functions and the theory of computation,
in \cite{La10}, 105--128.

\bibitem{WN09} D.\ Woods and T.\ Neary, The complexity of small universal Turing machines:
a survey, \emph{Theoret. Comput. Sci.} {\bf 410}, 2009, 443--450.
Extended and updated in: http://arxiv/abs/1110.2230.


\end{thebibliography}
\end{document}